\documentclass{article}

\usepackage{amsmath,amssymb,amsfonts,amsthm,mathtools,bm}
\usepackage[a4paper,left=2.8cm,right=2.8cm,top=2.6cm,bottom=2.7cm]{geometry}
\usepackage{microtype}
\usepackage{graphicx}
\usepackage{float}
\usepackage{needspace,placeins}
\usepackage{xcolor}
\usepackage{hyperref}

\definecolor{linkgreen}{RGB}{10,74,38}
\hypersetup{colorlinks=true,linkcolor=linkgreen,citecolor=linkgreen,urlcolor=linkgreen}
\numberwithin{equation}{section}

\theoremstyle{plain}
\newtheorem{theorem}{Theorem}[section]
\newtheorem{corollary}[theorem]{Corollary}

\newtheorem{lemma}[theorem]{Lemma}
\theoremstyle{definition}

\theoremstyle{remark}

\newenvironment{keywords}{{\noindent\bfseries Keywords.}}{\par}
\newenvironment{AMS}{{\noindent\bfseries 2020 Mathematics Subject Classification.}}{\par}

\newcommand{\Sd}{\mathbb S^d}
\newcommand{\R}{\mathbb R}
\newcommand{\dd}{\,\mathrm d}
\newcommand{\norm}[1]{\left\lVert #1\right\rVert}

\newcommand{\wce}{\operatorname{wce}}
\newcommand{\dist}{\operatorname{dist}}
\newcommand{\e}{\mathrm e}
\newcommand{\cH}{\mathcal H}
\newcommand{\cP}{\mathcal P}
\newcommand{\cX}{\mathcal X}
\newcommand{\cY}{\mathcal Y}

\newcommand{\PP}{\mathbb P}
\newcommand{\OpQ}{\mathcal Q}
\newcommand{\OpI}{\mathcal I}

\begin{document}

\title{Perturbed QMC designs and reduced QMC strength}
\author{Hao-Ning Wu\thanks{School of Mathematical Sciences and Fujian Key Laboratory of Mathematical Modeling and High Performance Scientific Computing, Xiamen University, Xiamen, Fujian 361005, China (\texttt{hnwu@xmu.edu.cn})}}
\date{}
\maketitle

\begin{abstract}
Small errors in sampling locations can change the convergence rate of a cubature rule. We study this effect for equal-weight integration on spheres and determine how accurately the nodes must be realized to retain optimal Sobolev rates. Under an average energy condition on the reference nodes, we establish sharp perturbation bounds and characterize how the decay of node errors affects quasi-Monte Carlo (QMC) strength, which measures the smoothness range over which optimal rates hold. We identify when this strength is preserved and when it can decrease, and show that sufficiently slow decay of node errors can destroy all optimal Sobolev rates. For a reference rule with optimal convergence, we further establish a necessary and sufficient condition on node accuracy for retaining the optimal rate uniformly over all permitted perturbations. Numerical experiments with spherical designs, Fekete points, and Coulomb-energy points illustrate the different stability regimes.
\end{abstract} 
 


\section{Introduction}

Although a cubature rule prescribes its sampling nodes, their locations may be perturbed in practice. These perturbations can limit the improvement in integration accuracy as the number of nodes increases. In this work, we study how accurately the nodes must be realized to retain optimal convergence rates for equal-weight cubature on spheres.

Let $\Sd:=\{x\in\mathbb{R}^{d+1}:\|x\|_2=1\}\subset\R^{d+1}$ for $d\ge2$ be the unit $d$-sphere, with normalized surface measure $\sigma_d$. Let 
\[\cX_N=\{x_{N,1},\ldots,x_{N,N}\}\]
be a set of $N$ distinct nodes on $\Sd$. For $s>d/2$, a sequence $\cX:=\{\cX_N\}$ is called a \emph{quasi-Monte Carlo} (QMC) \emph{design} sequence \cite{BrauchartSaffSloanWomersley2014} for the Sobolev space $H^s(\Sd)$ if
\begin{equation}\label{eq:intro-qmc-rate}
 \sup_{\norm{f}_{H^s}\le1}\left|\frac1N\sum_{j=1}^N f(x_{N,j})-\int_{\Sd}f\,\dd\sigma_d\right|\le C_{s,d}N^{-s/d},
\end{equation}
where $N$ ranges over an infinite subset of the positive integers and $C_{s,d}$ is independent of $N$. The rate $N^{-s/d}$ is optimal, as established in \cite{BrauchartHesse2007,Hesse2006Lower,HesseSloan2005WCE}, and its attainment does not require polynomial exactness. For a sequence with the QMC property \eqref{eq:intro-qmc-rate} for at least one $s>d/2$, its \emph{QMC strength} $s_{\cX}^*$ is the supremum of the Sobolev indices for which \eqref{eq:intro-qmc-rate} holds.

We determine how much of the reference QMC strength can be guaranteed
uniformly over all node realizations satisfying a prescribed accuracy.
Specifically, for each such $N$, let
\[\cY_N=\{y_{N,1},\ldots,y_{N,N}\}\]
be a perturbed configuration of $N$ distinct nodes, paired so that
\begin{equation}\label{eq:intro-displacement}
 \max_{1\le j\le N}\dist(x_{N,j},y_{N,j})\le\delta_N,
\end{equation}
where $\delta_N\ge0$ is a prescribed upper bound on the paired node displacements and $\dist(x,y)=\arccos(x\cdot y)$ is the geodesic distance between $x,y\in\Sd$. Let $\cY:=\{\cY_N\}$ be another sequence of point sets, indexed by the same node counts as $\cX$. The cubature weights remain fixed at $1/N$. Whenever $\cY$ has the QMC property \eqref{eq:intro-qmc-rate}, we denote its strength by $s_{\cY}^*$.

Our geometric hypothesis is the average energy bound
\begin{equation}\label{eq:average-energy}
 \frac1{N^2}\sum_{i\ne j}
 \dist(x_{N,i},x_{N,j})^{-\gamma}\le M_\gamma
\end{equation}
for every $0<\gamma<2$ and every node count $N$ in the sequence, with $M_\gamma$ independent of $N$. Condition \eqref{eq:average-energy} allows individual nodes to be closer than $N^{-1/d}$. A convenient sufficient condition ensuring this bound \eqref{eq:average-energy} is that the reference sets in $\cX$ are well separated. Let 
\[q(\cX_N):=\min_{i\ne j}\dist(x_{N,i},x_{N,j})\] 
be the separation of $\cX_N$. We show that
\begin{equation}\label{equ:separation}
 q(\cX_N)\ge c_0N^{-1/d}
\end{equation}
with $c_0>0$ independent of $N$ implies \eqref{eq:average-energy}.

Under the hypothesis \eqref{eq:average-energy}, we establish that the change in the cubature rule satisfies
the sharp estimate
\begin{equation}\label{eq:intro-modulus}
\begin{aligned}
 \sup_{\max\limits_j\dist(x_{N,j},y_{N,j})\le\delta}&
 \sup_{\norm{f}_{H^s}\le1}
 \left|\frac1N\sum_{j=1}^N f(y_{N,j})
       -\frac1N\sum_{j=1}^N f(x_{N,j})\right|\\
 &\asymp
 \begin{cases}
  \delta+N^{-1/2}\delta^{s-d/2},&d/2<s<d/2+1,\\
  \delta+N^{-1/2}\delta\sqrt{1+|\ln\delta|},&s=d/2+1,\\
  \delta,&s>d/2+1,
 \end{cases}
 \end{aligned}
\end{equation}
The estimate holds for every fixed $s>d/2$ and $0<\delta<\pi$. The constants may depend on $d$, $s$, and the corresponding energy bound, but not on $N$ or
$\delta$. For $s>d/2+1$, the energy condition \eqref{eq:average-energy} is not required for establishing \eqref{eq:intro-modulus}. 

To describe how the perturbation decay determines the retained QMC strength, we specialize to displacement bounds in \eqref{eq:intro-displacement} of the form
\begin{equation}\label{equ:displacement}
\delta_N=O(N^{-\beta})
\end{equation}
for some $\beta>0$. Smaller $\beta$ permits larger perturbations. Under the hypothesis \eqref{eq:average-energy}, estimate \eqref{eq:intro-modulus}
implies that, when $\beta>1/2$, for every $d/2<s\le\beta d$, $\cY$ is a QMC design sequence for $H^s(\Sd)$ if and only if $\cX$ is. The following theorem gives the resulting QMC-strength bound and establishes its sharpness.

\begin{theorem}[Sharp perturbation law for QMC designs]\label{thm:main}
Let $\beta>0$ and let $\cX=\{\cX_N\}$ be a QMC design sequence on $\Sd$ with $s_{\cX}^*>d/2$.
\begin{enumerate} 
\item[(i)] Suppose that $\beta>1/2$ and \eqref{eq:average-energy} holds.
Every paired sequence $\cY$ satisfying
\[
 \max_{1\le j\le N}\dist(x_{N,j},y_{N,j})=O(N^{-\beta})
\]
obeys
\begin{equation*}
 s_{\cY}^*\ge\min\{s_{\cX}^*,\beta d\}.
\end{equation*} 
If $s_{\cX}^*<\beta d$, every such perturbed sequence $\cY$ satisfies $s_{\cY}^*=s_{\cX}^*$. The lower bound is sharp: if $s_{\cX}^*\ge\beta d$, there exists a paired sequence $\cY$ with maximal displacement $\asymp N^{-\beta}$
and $s_{\cY}^*=\beta d$.
\item[(ii)] If $0<\beta\le1/2$, there exists a paired sequence $\cY$ with
maximal displacement $\asymp N^{-\beta}$ for which the QMC property
\eqref{eq:intro-qmc-rate} fails for every $s>d/2$.
\end{enumerate}
\end{theorem}  

The perturbation estimates also give a node-accuracy criterion at each
fixed smoothness index. For $s>d/2$ and $\delta\ge0$, let
\[
 E_s(\cX_N,\delta)
 :=\sup_{\max_j\dist(x_{N,j},y_{N,j})\le\delta}
   \sup_{\norm{f}_{H^s}\le1}
   \left|\frac1N\sum_{j=1}^N f(y_{N,j})
         -\int_{\Sd}f\,\dd\sigma_d\right|,
\]
where the outer supremum is over paired configurations of $N$ distinct
nodes.

\begin{corollary}[Cubature error and node accuracy]\label{cor:accuracy}
Let $\{\cX_N\}$ be a QMC design sequence for a fixed $H^s(\Sd)$, $s>d/2$.
Suppose that \eqref{eq:average-energy} holds if $s\le d/2+1$.
Then
\begin{equation}\label{eq:accuracy-balance}
 E_s(\cX_N,\delta)\asymp N^{-s/d}+\delta,\qquad0\le\delta<\pi,
\end{equation}
with constants independent of $N$ and $\delta$. For any $\delta_N\ge0$,
\begin{equation}\label{eq:accuracy-iff}
 E_s(\cX_N,\delta_N)=O(N^{-s/d})
 \quad\text{if and only if}\quad\delta_N=O(N^{-s/d}).
\end{equation}
\end{corollary}

At fixed node accuracy $0<\varepsilon<\pi$, the two terms in
\eqref{eq:accuracy-balance} balance at $N\asymp\varepsilon^{-d/s}$;
beyond this scale the worst error is of order $\varepsilon$.

The strength law also applies to \emph{spherical $t$-designs} \cite{DelsarteGoethalsSeidel1977}. A finite set $\cX_t\subset\Sd$ with $N_t:=|\cX_t|$ is said to be a \emph{spherical $t$-design} if
\begin{equation}\label{eq:intro-design}
 \frac1{N_t}\sum_{x\in\cX_t}p(x)
 =\int_{\Sd}p\,\dd\sigma_d\qquad \forall p\in \mathbb{P}_t,
\end{equation}
where $\mathbb{P}_t$ is the space of spherical polynomials of degree at most $t$. The Sobolev cubature bound established in \cite{BrauchartHesse2007} gives, for every fixed $s>d/2$,
\begin{equation}\label{eq:intro-design-error}
 \sup_{\norm{f}_{H^s}\le1} \left|\frac1{N_t}\sum_{x\in\cX_t}f(x)-\int_{\Sd}f\,\dd\sigma_d\right|\le C_{s,d}t^{-s},
\end{equation}
with $C_{s,d}$ independent of $t$. If $N_t\asymp t^d$, then $t^{-s}\asymp N_t^{-s/d}$, so the same family has the QMC property \eqref{eq:intro-qmc-rate} for every $s>d/2$, and its strength is infinite.
Algebraically decaying node errors can reduce this infinite strength to
a finite value.
The existence of well-separated spherical $t$-design families with $N_t\asymp t^d$ has been established in \cite{BondarenkoRadchenkoViazovska2013,BondarenkoRadchenkoViazovska2015}.
For these families, we express the maximal paired displacement as
$O(t^{-\alpha})$, where $\alpha>0$ is the decay exponent with respect
to $t$. Since $N_t\asymp t^d$, this corresponds to $\beta=\alpha/d$
in Theorem~\ref{thm:main}, which gives the following consequence.

\begin{corollary}[Perturbed spherical designs]\label{cor:designs}
Let $\alpha>0$ and let $\{\cX_t\}$ be a family of spherical $t$-designs
with $N_t\asymp t^d$.
\begin{enumerate}
\item[(i)] Suppose that $\alpha>d/2$ and \eqref{eq:average-energy} holds
with $N=N_t$. Every paired sequence $\cY=\{\cY_t\}$ with maximal
displacement $O(t^{-\alpha})$ is a QMC design sequence for $H^s(\Sd)$
for every $d/2<s\le\alpha$, and satisfies $s_{\cY}^*\ge\alpha$.
The lower bound is sharp: there exists a paired sequence $\cY$ with
maximal displacement $\asymp t^{-\alpha}$ and $s_{\cY}^*=\alpha$.
If the reference family is well separated, this sequence can also be
chosen well separated.
\item[(ii)] If $0<\alpha\le d/2$, there exists a paired sequence $\cY$
with maximal displacement $\asymp t^{-\alpha}$ for which the QMC property
\eqref{eq:intro-qmc-rate} fails for every $s>d/2$.
\end{enumerate}
\end{corollary}

\begin{proof}
Apply Theorem~\ref{thm:main} and the QMC equivalence preceding it
with $\beta=\alpha/d$.
If $q(\cX_t)\ge c_0t^{-1}$, a paired displacement of order
$t^{-\alpha}$ with $\alpha>d/2\ge1$ also gives
\[
 q(\cY_t)\ge q(\cX_t)-Ct^{-\alpha}
 \ge\frac{c_0}{2}t^{-1}
\]
for all sufficiently large $t$, proving the separation assertion.
\end{proof}

Optimal Sobolev cubature bounds based on polynomial exactness have been established in \cite{BrauchartHesse2007,HesseSloan2005WCE,HesseSloan2006ArbitraryOrder}, with matching complexity bounds for Besov classes in \cite{HesseMhaskarSloan2007}. The kernel formulation of QMC designs in
\cite{BrauchartSaffSloanWomersley2014} connects these errors with energy and discrepancy; see also \cite{BilykDaiMatzke2018,BrauchartDick2013Sobolev,BrauchartGrabner2015}. The relation between Sobolev cubature errors and covering radii has been studied in \cite{BrauchartEtAl2015Covering}. Beyond polynomial exactness, QMC designs have been characterized by finite harmonic sums in \cite{SloanWomersley2026}, nested equal-weight sequences have been constructed for each fixed $s>d/2$ in \cite{WuZhuang2026NestedQMC}. For random configurations, average QMC rates for stratified sampling have been established in \cite{BrauchartSaffSloanWomersley2014}, while average strengths of several point processes have been determined in \cite{DeLaTorreMarzo2024}. QMC rules also support approximation from samples, including hyperinterpolation, spectral multipliers, and spherical needlets \cite{AnWu2024Bypassing,BrauchartGrabnerSloanWomersley2024Needlets,Wu2026Hyperinterpolation}. Node perturbations have also been used in \cite{BrandoliniEtAl2014Quadrature} to demonstrate that a cubature error bound at one smoothness index need not extend to higher smoothness. In contrast, we quantify the sharp worst-case effect of uniformly bounded paired node displacements on equal-weight cubature and determine the resulting QMC-strength transition. Perturbed nodes also occur in spherical $t_\varepsilon$-designs and neighborhood cubature rules considered in \cite{Zhou2020Neighborhood,ZhouChen2018}, where positive weights are adjusted to preserve polynomial exactness. Here the equal weights remain fixed, and the question is which optimal Sobolev rates survive the node perturbations.

The paper is organized as follows. Section~\ref{sec:setup} collects the preliminaries. Section~\ref{sec:qmc} establishes the sharp perturbation estimates and proves Theorem~\ref{thm:main} and Corollary~\ref{cor:accuracy}. Section~\ref{sec:numerics} presents the numerical experiment.

\section{Preliminaries}\label{sec:setup}

Throughout, $d\ge2$ is fixed. $A\lesssim B$, or $B\gtrsim A$, means $A\le CB$ for an implicit constant $C>0$, and $A\asymp B$ means that both inequalities hold. Constants in $\lesssim$ and $\asymp$ may depend on $d$, fixed smoothness indices, and the QMC or energy bounds in the hypotheses, but not on $N$ or the perturbation radius. We let $x^{(k)}$ be the $k$th coordinate of $x$ and $p_k(x)=x^{(k)}$ for the
corresponding coordinate function.

For $\ell\ge0$, let $\cH_\ell(\Sd)$ be the space of restrictions to $\Sd$ of real homogeneous harmonic polynomials of degree $\ell$ on $\R^{d+1}$,
and set $Z(d,\ell):=\dim\cH_\ell(\Sd)$. Choose an orthonormal basis $\{Y_{\ell,k}:1\le k\le Z(d,\ell)\}$ of $\cH_\ell(\Sd)$ in $L_2(\Sd;\sigma_d)$.
Each $Y_{\ell,k}$ is an eigenfunction of the negative Laplace--Beltrami operator $-\Delta_{\Sd}$ with eigenvalue $\lambda_\ell=\ell(\ell+d-1)$.
The space
\[\PP_n(\Sd)=\bigoplus_{\ell=0}^n\cH_\ell(\Sd)\] 
is the space of spherical polynomials of degree at most $n$. For background, we refer the reader to \cite{AtkinsonHan2012,DaiXu2013}.
We define the Sobolev norm
\begin{equation}\label{eq:sobolev-norm}
 \norm{f}_{H^s}^2=\sum_{\ell=0}^\infty(1+\lambda_\ell)^s\sum_{k=1}^{Z(d,\ell)}|\widehat f_{\ell,k}|^2,
\end{equation}
where 
\[\widehat f_{\ell,k}=\langle f,Y_{\ell,k}\rangle = \int_{\Sd}fY_{\ell,k}d\sigma_d.\]

For a spherical integral
\[\OpI(f)=\int_{\Sd}f\,\dd\sigma_d\]
and an $N$-point set $\cX=\{x_1,\ldots,x_N\}$, the equal-weight cubature rule for the integral is given by
\begin{equation}\label{equ:cubature}
 \OpQ_{\cX}(f)=\frac1N\sum_{j=1}^N f(x_j).
\end{equation}
For $s>d/2$, the worst-case cubature error is
\begin{equation}\label{eq:wce}
 \wce(\OpQ_{\cX};H^s):=\sup_{\norm{f}_{H^s}\le1}|\OpQ_{\cX}(f)-\OpI(f)|=\sup_{\norm{f}_{H^s}\le1}
 \left|\frac1N\sum_{j=1}^N f(x_j)-\int_{\Sd}f\,\dd\sigma_d\right|.
\end{equation}
If $\cX$ is a QMC design sequence for $H^s(\Sd)$, then it is also a QMC design sequence for $H^{s'}(\Sd)$ for every $d/2<s'\le s$; see \cite[Theorem~4]{BrauchartSaffSloanWomersley2014}.

The following lemma justifies our geometric hypothesis \eqref{eq:average-energy}.
\begin{lemma}\label{lem:riesz}
Let $N\ge2$. If $q(\cX_N)\ge c_0N^{-1/d}$ and $0<\gamma<d$, then
\begin{equation}\label{eq:riesz-sum}
 \frac1{N^2}\sum_{i\ne j}\dist(x_i,x_j)^{-\gamma}
 \le C_{d,\gamma,c_0}.
\end{equation}
In particular, well-separated sequences satisfy
\eqref{eq:average-energy}.
\end{lemma}

\begin{proof}
Denoting by $q=q(\cX_N)$, we choose $K$ so that $2^Kq\le\pi<2^{K+1}q$. For fixed $i$ and $0\le k\le K$, the geodesic balls $B(x_j,q/2)$ with $2^kq\le\dist(x_i,x_j)<2^{k+1}q$ are disjoint and contained in $B(x_i,2^{k+1}q+q/2)$. Since each has area comparable to $q^d$,
\[
 \bigl|\{x_j:2^kq\le\dist(x_i,x_j)<2^{k+1}q\}\bigr|\lesssim\frac{\sigma_d(B(x_i,2^{k+1}q+q/2))}{q^d}\lesssim\frac{(2^kq)^d}{q^d}=2^{kd}.
\]
Grouping the points by these annuli gives
\[
 \begin{aligned}
 \sum_{j\ne i}\dist(x_i,x_j)^{-\gamma}
 &=\sum_{k=0}^K
   \sum_{\substack{j\ne i\\2^kq\le\dist(x_i,x_j)<2^{k+1}q}}
   \dist(x_i,x_j)^{-\gamma}\\
 &\le\sum_{k=0}^K(2^kq)^{-\gamma}
   \bigl|\{x_j:2^kq\le\dist(x_i,x_j)<2^{k+1}q\}\bigr|\\
 &\lesssim q^{-\gamma}\sum_{k=0}^K2^{k(d-\gamma)}
 \lesssim q^{-d}\le c_0^{-d}N.
 \end{aligned}
\]
Sum over $i$ and divide by $N^2$.
\end{proof}

\section{Proof of the QMC-strength law}\label{sec:qmc}\label{sec:sobolev}

For $\cX_N=\{x_1,\ldots,x_N\}$, we define
\[
 \cP_\delta(\cX_N)
 =\{\cY=\{y_1,\ldots,y_N\}\subset\Sd:
       \max_j\dist(x_j,y_j)\le\delta\}.
\]
The $y_j$ are distinct and retain the weights $1/N$.
For $s>d/2$, the perturbation error is
\begin{equation}\label{eq:perturbation-modulus}
 \omega_s(\cX_N,\delta)
 :=\sup_{\cY\in\cP_\delta(\cX_N)}\sup_{\norm{f}_{H^s}\le1}
       \left|\frac1N\sum_{j=1}^N(f(y_j)-f(x_j))\right|.
\end{equation}
For fixed $s>d/2$, the Sobolev norm \eqref{eq:sobolev-norm} has the
reproducing kernel
\[
 K(x,y)=\sum_{\ell=0}^{\infty}(1+\lambda_\ell)^{-s}
        \sum_{k=1}^{Z(d,\ell)}Y_{\ell,k}(x)Y_{\ell,k}(y),
 \qquad x,y\in\Sd;
\]
see \cite[Section~2.4, (17)--(18)]{BrauchartSaffSloanWomersley2014}.
Below, we also use kernels inducing equivalent Sobolev norms.
For such a kernel $K$, write $\langle\cdot,\cdot\rangle_K$ and
$\norm{\cdot}_K$ for the associated inner product and norm. Then
\[
 K(\cdot,x)\in H^s(\Sd),\qquad
 f(x)=\langle f,K(\cdot,x)\rangle_K,\qquad
 \norm{f}_K\asymp\norm{f}_{H^s}.
\]
For $\cY\in\cP_\delta(\cX_N)$, it follows that
\[
 \OpQ_{\cY}(f)-\OpQ_{\cX_N}(f)
 =\left\langle f,\frac1N\sum_{j=1}^N
       \bigl(K(\cdot,y_j)-K(\cdot,x_j)\bigr)\right\rangle_K.
\]
Cauchy--Schwarz, with equality for a normalized kernel sum when it is nonzero, gives the supremum over the $\norm{\cdot}_K$ unit ball.
Expanding its square and using $\langle K(\cdot,u),K(\cdot,v)\rangle_K=K(u,v)$ yields
\begin{equation}\label{eq:energy-expanded}
 \begin{aligned}
 \sup_{\norm{f}_{H^s}\le1}
       |\OpQ_{\cY}(f)-\OpQ_{\cX_N}(f)|^2&\asymp\sup_{\norm{f}_K\le1}
       |\OpQ_{\cY}(f)-\OpQ_{\cX_N}(f)|^2\\
 &=\left\|\frac1N\sum_{j=1}^N
       \bigl(K(\cdot,y_j)-K(\cdot,x_j)\bigr)\right\|_K^2\\
 &=\frac1{N^2}\sum_{i,j=1}^N
 \underbrace{[K(y_i,y_j)-K(y_i,x_j)-K(x_i,y_j)+K(x_i,x_j)]}_{=:B_{ij}}.
 \end{aligned}
\end{equation}
This is the argument in \cite[Section~3.2]{BrauchartSaffSloanWomersley2014}, with the integral replaced by the reference rule $\OpQ_{\cX_N}$.
We let $s_c=d/2+1$ be the critical index.

\subsection{Subcritical case}\label{sec:negative}

\begin{lemma}\label{lem:coordinate-perturbation}
For every $s>d/2$ and every $N$-point set $\cX_N$,
\begin{equation}\label{eq:modulus-linear-lower}
 \omega_s(\cX_N,\delta)\ge c_{d,s}\delta,
 \qquad 0<\delta<\pi.
\end{equation}
\end{lemma}

\begin{proof}
First suppose $0<\delta\le\delta_0$, where $0<\delta_0<1/2$ will be
chosen depending only on $d$. We construct a perturbation that changes
the average of one coordinate by an amount comparable to $\delta$.
Let $p_k(x)=x^{(k)}$ be the $k$th coordinate of $x\in\mathbb R^{d+1}$.
Since $\sum_{k=1}^{d+1}(1-(x^{(k)})^2)=d$, we can choose $k$ such that
\begin{equation}\label{eq:coordinate-gradient-average}
 \OpQ_{\cX_N}(1-p_k^2)
 =\frac1N\sum_{j=1}^N(1-(x_j^{(k)})^2)\ge\frac d{d+1}.
\end{equation}
To perturb the nodes, we translate in the $k$th coordinate direction
and normalize back to the sphere. For the $k$th coordinate vector $e_k$,
we let, for $|\eta|<1/2$,
\begin{equation}\label{eq:coordinate-map}
 T_{k,\eta}(x):=\frac{x+\sin\eta\,e_k}{|x+\sin\eta\,e_k|}.
\end{equation}
The origin lies inside the ball bounded by the translated sphere
$\Sd+\sin\eta\,e_k$. Each ray from the origin meets this sphere
exactly once, so radial projection, and hence $T_{k,\eta}$, is injective.
To bound the displacement, we denote by $a=\sin\eta$, $u=x^{(k)}$,
and $\varphi=\dist(x,T_{k,\eta}(x))$. Then $0\le\varphi<\pi/2$ and
\[
 \sin^2\varphi=\frac{a^2(1-u^2)}{1+2au+a^2}\le a^2,
\]
where the equality follows from $\sin^2\varphi+\cos^2\varphi=1$
and the inequality follows from $(u+a)^2\ge0$.
Thus $\varphi\le\arcsin|a|=|\eta|$.
Using $|a|\asymp|\eta|$ and $\sqrt{1+2au+a^2}\le1+|a|$ in the
same formula also gives
\begin{equation}\label{eq:coordinate-displacement}
 c|\eta|\sqrt{1-(x^{(k)})^2}
 \le\dist(x,T_{k,\eta}(x))\le|\eta|,
\end{equation}
with an absolute constant $c>0$.
Expanding the $k$th coordinate of $T_{k,\eta}(x)$ at $\eta=0$ gives
\begin{equation}\label{eq:coordinate-taylor}
 p_k(T_{k,\eta}(x))-p_k(x)
 =\eta(1-(x^{(k)})^2)+O(\eta^2),
\end{equation}
uniformly in $x$ and $k$. These properties hold for every coordinate
index $k$ and both signs of $\eta$.
Now take $\cY=T_{k,\delta}(\cX_N)$, with $k$ chosen in
\eqref{eq:coordinate-gradient-average}. Injectivity and
\eqref{eq:coordinate-displacement} give $\cY\in\cP_\delta(\cX_N)$.
Averaging \eqref{eq:coordinate-taylor} with $\eta=\delta$ and using
\eqref{eq:coordinate-gradient-average} gives
\[
 (\OpQ_{\cY}-\OpQ_{\cX_N})(p_k)
 =\delta\OpQ_{\cX_N}(1-p_k^2)+O(\delta^2)\ge\frac d{d+1}\delta-C\delta^2\ge c_d\delta
\]
for a sufficiently small $\delta_0$ depending only on $d$.
Dividing by $\norm{p_k}_{H^s}$ and using
\eqref{eq:perturbation-modulus} proves \eqref{eq:modulus-linear-lower}
for $0<\delta\le\delta_0$.
For $\delta_0<\delta<\pi$, the inclusion
$\cP_{\delta_0}(\cX_N)\subset\cP_\delta(\cX_N)$ gives
\[
 \omega_s(\cX_N,\delta)\ge\omega_s(\cX_N,\delta_0)
 \ge c_{d,s}\delta_0\ge\frac{c_{d,s}\delta_0}{\pi}\,\delta.
\]
Since $\delta_0$ depends only on $d$, adjusting $c_{d,s}$ proves the
claimed bound for every $0<\delta<\pi$.
\end{proof}

For $d/2<s<s_c$, put $\theta=2s-d\in(0,2)$. An equivalent
reproducing kernel of $H^s(\Sd)$ is
\begin{equation}\label{eq:gd-kernel}
 K_s(x,y)=C_{s,d}-|x-y|^\theta,
\end{equation}
with $C_{s,d}$ sufficiently large;
see \cite[Section~5]{BrauchartSaffSloanWomersley2014}.
We use $|x-y|=2\sin(\dist(x,y)/2)\asymp\dist(x,y)$.

\begin{lemma}\label{lem:mixed}
Let $x_i,x_j,y_i,y_j\in\Sd$, $0<\theta<2$, and
$\dist(x_i,x_j)>0$. If 
\[\max\{\dist(x_i,y_i),\dist(x_j,y_j)\}\le\delta,\] 
then
\begin{equation}\label{eq:mixed-difference}
 \bigl||y_i-y_j|^\theta-|y_i-x_j|^\theta
       -|x_i-y_j|^\theta+|x_i-x_j|^\theta\bigr|
 \le C_{d,\theta}\delta^2\dist(x_i,x_j)^{\theta-2}.
\end{equation}
\end{lemma}

\begin{proof}
Let $z=x_i-x_j$, $a=y_i-x_i$, $b=y_j-x_j$, and $\Phi(w)=|w|^\theta$. If $\dist(x_i,x_j)\ge4\delta$, then for $0\le u,v\le1$,
\[\begin{aligned}
 |z+ua-vb|
 &\ge |z|-u|a|-v|b|\ge |x_i-x_j|-2\delta\\
 &\ge\frac{2\dist(x_i,x_j)}{\pi}-2\delta\\
 &\ge\left(\frac2\pi-\frac12\right)\dist(x_i,x_j).
 \end{aligned}
\]
Here the third inequality uses $|x_i-x_j|=2\sin(\dist(x_i,x_j)/2)$
and $\sin t\ge2t/\pi$ for $0\le t\le\pi/2$, while the last uses
$2\delta\le\dist(x_i,x_j)/2$.
The mixed difference is
\[
 \Phi(z+a-b)-\Phi(z+a)-\Phi(z-b)+\Phi(z)
 =-\int_0^1\int_0^1a^\top D^2\Phi(z+ua-vb)b\,\dd u\dd v.
\]
Here $D^2\Phi$ is the Euclidean Hessian, with operator norm
$\norm{D^2\Phi(w)}\le C_\theta|w|^{\theta-2}$ for $w\ne0$.
This gives the bound. If $\dist(x_i,x_j)<4\delta$, the triangle inequality gives
\[
 \begin{aligned}
 |y_i-y_j|\le |y_i-x_i|+|x_i-x_j|+|x_j-y_j|\le\delta+\dist(x_i,x_j)+\delta<6\delta.
 \end{aligned}
\]
The other three distances are at most $5\delta$, $5\delta$, and
$4\delta$, respectively. Hence the absolute value of the mixed difference
is at most $4(6\delta)^\theta=C_\theta\delta^\theta$.
Since $2-\theta>0$,
\[
 \delta^\theta=\delta^2\dist(x_i,x_j)^{\theta-2}
   \left(\frac{\dist(x_i,x_j)}{\delta}\right)^{2-\theta}\le4^{2-\theta}\delta^2\dist(x_i,x_j)^{\theta-2},
\]
which proves the same bound in this case.
\end{proof}

\begin{theorem}[Subcritical perturbation estimate]\label{thm:negative}
Let $d/2<s<s_c$ and suppose that
\begin{equation}\label{eq:average-s-energy}
 \frac1{N^2}\sum_{i\ne j}\dist(x_i,x_j)^{-(d+2-2s)}\le M_s.
\end{equation}
Then, for $0<\delta<\pi$,
\begin{equation}\label{eq:negative-stability}
 \omega_s(\cX_N,\delta)
 \asymp\delta+N^{-1/2}\delta^{s-d/2}.
\end{equation}
The comparison constants depend only on $d$, $s$, and $M_s$.
\end{theorem}

\begin{proof}
\underline{\emph{Upper bound.}} For $\cY\in\cP_\delta(\cX_N)$ and $K=K_s$, the summands in
\eqref{eq:energy-expanded} are
\[
 B_{ij}=|y_i-x_j|^\theta+|x_i-y_j|^\theta
        -|y_i-y_j|^\theta-|x_i-x_j|^\theta,
 \qquad B_{ii}=2|x_i-y_i|^\theta.
\]
Since $|x_i-y_i|\le\dist(x_i,y_i)\le\delta$, the diagonal contribution satisfies
\[
 \frac1{N^2}\sum_{i=1}^N B_{ii}
 =\frac2{N^2}\sum_{i=1}^N|x_i-y_i|^\theta
 \le\frac{2N\delta^\theta}{N^2}=2N^{-1}\delta^\theta.
\]
For $i\ne j$, apply Lemma~\ref{lem:mixed} and note that
$\theta-2=-(d+2-2s)$. The energy bound \eqref{eq:average-s-energy} then gives
\[
 \frac1{N^2}\sum_{i\ne j}|B_{ij}|\le\frac{C_{d,s}\delta^2}{N^2}
       \sum_{i\ne j}\dist(x_i,x_j)^{-(d+2-2s)}\le C_{d,s}M_s\delta^2.
\]
Together with \eqref{eq:energy-expanded}, these bounds give
$\omega_s(\cX_N,\delta)^2
\lesssim N^{-1}\delta^\theta+\delta^2$.

\underline{\emph{Lower bound.}} Choose distinct $y_j$ with $\dist(x_j,y_j)=\delta$; this is possible
successively because $d\ge2$. Then
$|x_j-y_j|=2\sin(\delta/2)\ge2\delta/\pi$, so
\begin{equation}\label{eq:energy-lower-two-regimes}
 \frac1{N^2}\sum_{i,j=1}^N B_{ij}
 =\frac1{N^2}\sum_iB_{ii}+\frac1{N^2}\sum_{i\ne j}B_{ij}
 \ge cN^{-1}\delta^\theta-C\delta^2.
\end{equation}
For $\delta^{2-\theta}\le c'/N$ with $c'$ sufficiently small, this yields
$\omega_s(\cX_N,\delta)\gtrsim N^{-1/2}\delta^{\theta/2}$.
Otherwise $N^{-1/2}\delta^{\theta/2}\lesssim\delta$.
Combine these two cases with Lemma~\ref{lem:coordinate-perturbation}.
\end{proof}

\subsection{Critical case and higher smoothness}

At $s=s_c$, the power kernel \eqref{eq:gd-kernel} reduces to a
polynomial of degree one in $x\cdot y$. Taking a limit of the expansion in
\cite[Section~5, (25)--(26)]{BrauchartSaffSloanWomersley2014}
gives the following logarithmic kernel.

\begin{lemma}[A kernel at the critical index]\label{lem:critical-kernel}
With the convention $0^2\ln0:=0$ and $A_d,B_d>0$ sufficiently large,
\begin{equation}\label{eq:critical-kernel}
 K_c(x,y)=A_d+B_d\,x\cdot y+|x-y|^2\ln|x-y|
\end{equation}
is a reproducing kernel for a norm equivalent to that of $H^{s_c}(\Sd)$.
Define its spherical harmonic coefficients $b_\ell$ by
\[
 K_c(x,y)=\sum_{\ell=0}^{\infty}b_\ell
          \sum_{k=1}^{Z(d,\ell)}Y_{\ell,k}(x)Y_{\ell,k}(y).
\]
For $\ell\ge2$, these coefficients satisfy
\begin{equation}\label{eq:critical-coefficients}
 b_\ell=\frac{\Gamma(d+1)\Gamma(\ell-1)}{\Gamma(\ell+d+1)}
 \asymp(1+\lambda_\ell)^{-s_c},\qquad\ell\ge2.
\end{equation}
\end{lemma}

\begin{proof}
For $0<\theta<2$, let
$V_\theta=\iint|x-y|^\theta\,\dd\sigma_d(x)\dd\sigma_d(y)$.
The spherical harmonic expansion of $-|x-y|^\theta$ is
\[
 -|x-y|^\theta
 =-V_\theta-V_\theta\sum_{\ell=1}^{\infty}
   \frac{(-\theta/2)_\ell}{(d+\theta/2)_\ell}
   \sum_{k=1}^{Z(d,\ell)}Y_{\ell,k}(x)Y_{\ell,k}(y);
\]
see \cite[Section~5, (25)--(26)]{BrauchartSaffSloanWomersley2014}, where
$(a)_0=1$ and $(a)_\ell=a(a+1)\cdots(a+\ell-1)$ for $\ell\ge1$.
The term $A_d+B_d\,x\cdot y$ affects only degrees zero and one.
Thus, for $\ell\ge2$, $b_\ell$ is also the coefficient of
$\sum_{k=1}^{Z(d,\ell)}Y_{\ell,k}(x)Y_{\ell,k}(y)$ in the expansion
of $|x-y|^2\ln|x-y|$. We use the limit
\[
 \frac{|x-y|^2-|x-y|^{2-\varepsilon}}{\varepsilon}
 \longrightarrow |x-y|^2\ln|x-y|
 \quad(\varepsilon\rightarrow0^+),
\]
uniformly for $x,y\in\Sd$.
Uniform convergence allows the limit to pass through each fixed spherical
harmonic projection.
Since $|x-y|^2=2-2x\cdot y$ has no spherical harmonic components of
degree $\ell\ge2$, taking the coefficient in each fixed degree in this
limit and using the preceding expansion gives
\[
 b_\ell=\lim_{\varepsilon\rightarrow0^+}
 \frac{-V_{2-\varepsilon}(-1+\varepsilon/2)_\ell}
      {\varepsilon(d+1-\varepsilon/2)_\ell}
 =\frac{(\ell-2)!}{(d+1)_\ell}.
\]
Here $V_2=2$ and, for $\ell\ge2$,
$(-1+\varepsilon/2)_\ell
=-\varepsilon(\ell-2)!/2+O_\ell(\varepsilon^2)$.
All $b_\ell$ with $\ell\ge2$ are positive and decay as $\ell^{-d-2}$.
Adding $A_d+B_d\,x\cdot y$ makes the two remaining coefficients
positive. Since $\sum_\ell b_\ell Z(d,\ell)<\infty$, the expansion
converges absolutely and uniformly and defines the required Sobolev norm.
\end{proof}

\begin{theorem}[Critical and higher-smoothness estimates]\label{thm:critical-smooth}
\leavevmode
\begin{enumerate}
\item[(i)] If
\begin{equation}\label{eq:logarithmic-packing}
 \frac1{N^2}\sum_{i\ne j}(1+|\ln \dist(x_i,x_j)|)\le M_c,
\end{equation}
then, for $0<\delta<\pi$,
\begin{equation}\label{eq:critical-modulus}
 \omega_{s_c}(\cX_N,\delta)
 \asymp\delta+N^{-1/2}\delta\sqrt{1+|\ln\delta|}.
\end{equation}
The comparison constants depend only on $d$ and $M_c$.
\item[(ii)] For every $s>s_c$ and every $N$-point set $\cX_N$,
\begin{equation}\label{eq:smooth-modulus}
 \omega_s(\cX_N,\delta)\asymp\delta,\qquad0<\delta<\pi,
\end{equation}
with constants depending only on $d$ and $s$.
\end{enumerate}
\end{theorem}

\begin{proof}
\underline{\emph{Critical upper bound.}}
First suppose that $0<\delta\le1/12$ and use $K=K_c$ in \eqref{eq:energy-expanded}.
Since $K_c(x_i,x_i)=K_c(y_i,y_i)=A_d+B_d$ and
$|x_i-y_i|^2=2(1-x_i\cdot y_i)$, the diagonal terms satisfy
\begin{equation}\label{eq:critical-diagonal}
 \begin{aligned}
 B_{ii}
 &=K_c(y_i,y_i)-2K_c(x_i,y_i)+K_c(x_i,x_i)\\
 &=2(A_d+B_d)-2\bigl(A_d+B_d\,x_i\cdot y_i
                    +|x_i-y_i|^2\ln|x_i-y_i|\bigr)\\
 &=|x_i-y_i|^2\bigl(B_d-2\ln|x_i-y_i|\bigr).
 \end{aligned}
\end{equation}
Since $|x_i-y_i|\le\dist(x_i,y_i)\le\delta\le1/12$ and
\[
 \frac{\dd}{\dd t}\bigl(t^2\ln(1/t)\bigr)
 =t\bigl(2\ln(1/t)-1\bigr)>0,\qquad0<t\le1/2,
\]
we obtain, with the value at $x_i=y_i$ understood by continuity,
\[
 \begin{aligned}
 B_{ii}
 &=B_d|x_i-y_i|^2+2|x_i-y_i|^2\ln\frac1{|x_i-y_i|}\\
 &\le B_d\delta^2+2\delta^2\ln(1/\delta)
 \lesssim\delta^2(1-\ln\delta).
 \end{aligned}
\]
For $F(w)=|w|^2\ln|w|$ and $w\ne0$, differentiation gives
\[
 \nabla F(w)=(2\ln|w|+1)w,\qquad
 D^2F(w)v=(2\ln|w|+1)v+2\frac{w\cdot v}{|w|^2}w
\]
for every $v\in\mathbb R^{d+1}$. Hence, by $|w\cdot v|\le|w||v|$,
\[
 |D^2F(w)v|
 \le\bigl(|2\ln|w|+1|+2\bigr)|v|
 \le\bigl(3+2|\ln|w||\bigr)|v|.
\]
Taking the supremum over $|v|=1$ yields
$\norm{D^2F(w)}\lesssim1+|\ln|w||$.
If $\dist(x_i,x_j)\ge4\delta$, use $z,a,b$ from the proof of
Lemma~\ref{lem:mixed}. The path estimate there and the triangle inequality give
\[
 |z+ua-vb|\asymp\dist(x_i,x_j),\qquad0\le u,v\le1.
\]
Since $|a|,|b|\le\delta$, the integrand in the mixed-difference formula satisfies
\[
 |a^\top D^2F(z+ua-vb)b|
 \le\delta^2\norm{D^2F(z+ua-vb)}
 \lesssim\delta^2(1+|\ln\dist(x_i,x_j)|).
\]
Integration gives the same bound for the mixed difference of $F$.
If $\dist(x_i,x_j)<4\delta$, the four distances are at most
$6\delta\le1/2$, as in Lemma~\ref{lem:mixed}.
By the monotonicity of $t^2\ln(1/t)$ established above, each value satisfies
\[
 |F(w)|\le(6\delta)^2\ln\frac1{6\delta}
 \le36\delta^2\ln(1/\delta),\qquad |w|\le6\delta.
\]
Moreover, $\dist(x_i,x_j)<4\delta\le1/3$ gives
\[
 1-\ln\delta
 \le1+\ln4-\ln\dist(x_i,x_j)
 \lesssim1+|\ln\dist(x_i,x_j)|.
\]
Thus the mixed difference of $F$ has the same bound in this case.
The constant term $A_d$ cancels, while $B_d\,x\cdot y$ contributes
\[
 |B_d(y_i-x_i)\cdot(y_j-x_j)|\le B_d\delta^2.
\]
Combining these estimates yields
\begin{equation}\label{eq:critical-mixed}
 |B_{ij}|\lesssim\delta^2(1+|\ln \dist(x_i,x_j)|),\qquad i\ne j.
\end{equation}
The diagonal estimate above gives
\[
 \frac1{N^2}\sum_i B_{ii}
 \lesssim N^{-1}\delta^2(1-\ln\delta).
\]
For the off-diagonal terms, \eqref{eq:critical-mixed} and
\eqref{eq:logarithmic-packing} give
\[
 \frac1{N^2}\sum_{i\ne j}|B_{ij}|
 \lesssim\frac{\delta^2}{N^2}\sum_{i\ne j}
 (1+|\ln\dist(x_i,x_j)|)
 \le M_c\delta^2.
\]
Substituting into \eqref{eq:energy-expanded} and taking the supremum
over $\cY\in\cP_\delta(\cX_N)$ yields
\[
 \omega_{s_c}(\cX_N,\delta)^2
 \lesssim N^{-1}\delta^2(1-\ln\delta)+M_c\delta^2.
\]
For $1/12<\delta<\pi$, uniform boundedness of point evaluation gives
\[
 \omega_{s_c}(\cX_N,\delta)
 \le2\sup_{\norm{f}_{H^{s_c}}\le1}\norm{f}_\infty
 \le C_d\le12C_d\delta.
\]
Thus the upper bound in \eqref{eq:critical-modulus} holds for all
$0<\delta<\pi$.

\underline{\emph{Critical lower bound.}}
First suppose that $0<\delta\le1/12$. Choose distinct nodes at distance $\delta$ from their reference nodes,
as in Theorem~\ref{thm:negative}. Then $|x_i-y_i|=2\sin(\delta/2)$ and
\[
 B_{ii}=4\sin^2(\delta/2)
       \bigl(B_d-2\ln(2\sin(\delta/2))\bigr)
 \gtrsim\delta^2(1-\ln\delta).
\]
Together with \eqref{eq:critical-mixed} and
\eqref{eq:logarithmic-packing}, this gives
\[
 \frac1{N^2}\sum_{i,j=1}^N B_{ij}
 \ge cN^{-1}\delta^2(1-\ln\delta)-C\delta^2.
\]
If $N^{-1}(1-\ln\delta)\ge2C/c$, the right-hand side is at least
$(c/2)N^{-1}\delta^2(1-\ln\delta)$. Thus
\eqref{eq:energy-expanded} and \eqref{eq:perturbation-modulus} yield
\[
 \omega_{s_c}(\cX_N,\delta)
 \gtrsim N^{-1/2}\delta\sqrt{1-\ln\delta}.
\]
Combining this with the linear lower bound
\eqref{eq:modulus-linear-lower} controls the sum of the two terms in
\eqref{eq:critical-modulus}, since $\omega_{s_c}$ is a supremum over
all permitted realizations.
If $N^{-1}(1-\ln\delta)<2C/c$, then
\[
 N^{-1/2}\delta\sqrt{1-\ln\delta}
 \le\sqrt{2C/c}\,\delta,
\]
so \eqref{eq:modulus-linear-lower} alone controls both terms.
This proves the lower bound for $0<\delta\le1/12$, where
$1-\ln\delta=1+|\ln\delta|$.
For $1/12<\delta<\pi$, the same lemma gives
\[
 \omega_{s_c}(\cX_N,\delta)\gtrsim\delta.
\]
Since $N^{-1/2}\sqrt{1+|\ln\delta|}\lesssim1$ on this interval,
this proves the lower bound in \eqref{eq:critical-modulus} for all
$0<\delta<\pi$.

\underline{\emph{Higher smoothness.}}
For $s>s_c=d/2+1$, Sobolev embedding gives
\[
 |f(x)-f(y)|\le C_{d,s}\norm{f}_{H^s}\dist(x,y),
 \qquad x,y\in\Sd.
\]
Thus, for every $\cY\in\cP_\delta(\cX_N)$,
\[
 \begin{aligned}
 |(\OpQ_{\cY}-\OpQ_{\cX_N})(f)|
 &\le\frac1N\sum_{j=1}^N|f(y_j)-f(x_j)|\le C_{d,s}\norm{f}_{H^s}
       \frac1N\sum_{j=1}^N\dist(x_j,y_j)\\
 &\le C_{d,s}\delta\norm{f}_{H^s}.
 \end{aligned}
\]
Taking the suprema in \eqref{eq:perturbation-modulus} gives
$\omega_s(\cX_N,\delta)\lesssim\delta$.
The reverse inequality follows from \eqref{eq:modulus-linear-lower},
proving \eqref{eq:smooth-modulus}.
\end{proof}

Under \eqref{eq:average-energy}, the subcritical condition
\eqref{eq:average-s-energy} follows by taking $\gamma=d+2-2s$.
The critical condition follows from
$1+|\ln r|\lesssim1+r^{-1}$ on $(0,\pi]$ and $\gamma=1$.
Thus Theorems~\ref{thm:negative} and \ref{thm:critical-smooth}
give \eqref{eq:intro-modulus}; Lemma~\ref{lem:riesz} covers the
well-separated case.

\subsection{Proofs of Theorem~\texorpdfstring{\ref{thm:main}}{1.1} and Corollary~\texorpdfstring{\ref{cor:accuracy}}{1.2}}

\begin{proof}[Proof of Theorem~\ref{thm:main}]
\underline{\emph{Stability.}}
We first prove the strength lower bound in part~(i) and the equality
when $s_{\cX}^*<\beta d$. Let $\beta>1/2$ and let $\cY$ satisfy
\[\max_j\dist(x_{N,j},y_{N,j})\le C N^{-\beta}\] 
for all sufficiently large $N$, with $C$ independent of $N$. Fix $d/2<s\le\beta d$.
For $s<s_c$, Theorem~\ref{thm:negative} gives
\begin{equation}\label{eq:qmc-two-terms}
 \omega_s(\cX_N,CN^{-\beta})
 \lesssim N^{-\beta}+N^{-1/2-\beta(s-d/2)}
 \lesssim N^{-s/d},
\end{equation}
using $s/d\le\beta$ and
\begin{equation}\label{eq:exponent-identity}
 \frac12+\beta\left(s-\frac d2\right)-\frac sd
 =\left(\beta-\frac1d\right)\left(s-\frac d2\right)>0.
\end{equation}
At $s=s_c$, apply Theorem~\ref{thm:critical-smooth}(i) with
$\delta=C N^{-\beta}$ to obtain, for sufficiently large $N$,
\[
 \omega_{s_c}(\cX_N,CN^{-\beta})
 \lesssim N^{-\beta}\left(1+\sqrt{\frac{\ln N}{N}}\right)
 \lesssim N^{-s_c/d}.
\]
For $s>s_c$, the bound is directly $O(N^{-\beta})$.
For every $f$ with $\norm{f}_{H^s}\le1$, the triangle inequality gives
\[
 \begin{aligned}
 |\OpQ_{\cY_N}(f)-\OpI(f)|
 &\le |\OpQ_{\cX_N}(f)-\OpI(f)|
       +|\OpQ_{\cY_N}(f)-\OpQ_{\cX_N}(f)|\\
 &\le \wce(\OpQ_{\cX_N};H^s)
       +\omega_s(\cX_N,CN^{-\beta}),
 \end{aligned}
\]
where the last step uses \eqref{eq:perturbation-modulus} and
$\cY_N\in\cP_{CN^{-\beta}}(\cX_N)$.
Taking the supremum over $f$ yields the first of the following inequalities;
applying the same argument with the two rules interchanged yields the second:
\[
 \begin{aligned}
 \wce(\OpQ_{\cY_N};H^s)
 &\le\wce(\OpQ_{\cX_N};H^s)+\omega_s(\cX_N,CN^{-\beta}),\\
 \wce(\OpQ_{\cX_N};H^s)
 &\le\wce(\OpQ_{\cY_N};H^s)+\omega_s(\cX_N,CN^{-\beta}).
 \end{aligned}
\]
Rearranging these two inequalities gives
\[
 -\omega_s(\cX_N,CN^{-\beta})
 \le\wce(\OpQ_{\cY_N};H^s)-\wce(\OpQ_{\cX_N};H^s)
 \le\omega_s(\cX_N,CN^{-\beta}),
\]
yielding
\begin{equation}\label{eq:qmc-equivalence}
 \left|\wce(\OpQ_{\cY_N};H^s)-\wce(\OpQ_{\cX_N};H^s)\right|\le\omega_s(\cX_N,CN^{-\beta})\lesssim N^{-s/d},\qquad d/2<s\le\beta d.
\end{equation}
This proves the equivalence of the $H^s$-QMC properties of $\cX$ and
$\cY$ for every $d/2<s\le\beta d$, as stated immediately before
Theorem~\ref{thm:main}. It also gives
$s_{\cY}^*\ge\min\{s_{\cX}^*,\beta d\}$.
If $s_{\cX}^*<\beta d$, this lower bound gives
$s_{\cY}^*\ge s_{\cX}^*$. To prove equality, suppose for a contradiction
that $s_{\cY}^*>s_{\cX}^*$ and choose
\[
 s_{\cX}^*<s<\min\{s_{\cY}^*,\beta d\}.
\]
By the definition of strength and the fact that the QMC property passes
to lower smoothness indices, $s<s_{\cY}^*$ implies that $\cY$ is a QMC
design sequence for $H^s(\Sd)$. Since $s<\beta d$, the equivalence
\eqref{eq:qmc-equivalence} gives the same property for $\cX$.
This contradicts $s>s_{\cX}^*$, because $s_{\cX}^*$ is the supremum
of its admissible QMC indices. Hence $s_{\cY}^*\le s_{\cX}^*$, and
therefore $s_{\cY}^*=s_{\cX}^*$.

\underline{\emph{Sharpness and loss of QMC rates.}}
We now construct a perturbed sequence that attains the lower bound in
part~(i) when $s_{\cX}^*\ge\beta d$ and proves part~(ii).
For any $\beta>0$, choose $\varepsilon_N\in\{-1,1\}$ so that
$\varepsilon_N\OpQ_{\cX_N}(p_{d+1})
=|\OpQ_{\cX_N}(p_{d+1})|$. Using the map defined in
\eqref{eq:coordinate-map}, set, for all sufficiently large $N$,
\begin{equation}\label{eq:signed-perturbation}
 y_{N,j}=T_{d+1,\varepsilon_NN^{-\beta}}(x_{N,j}).
\end{equation}
Take $\cY_N=\cX_N$ at the remaining indices.
Since the reference sequence is QMC at some index greater than $d/2$,
\[
 \OpQ_{\cX_N}(1-p_{d+1}^2)\longrightarrow
 \OpI(1-p_{d+1}^2)=\frac d{d+1}.
\]
Equation~\eqref{eq:coordinate-taylor} gives
\begin{equation}\label{eq:signed-obstruction}
 \varepsilon_N\OpQ_{\cY_N}(p_{d+1})
 =|\OpQ_{\cX_N}(p_{d+1})|
  +N^{-\beta}\OpQ_{\cX_N}(1-p_{d+1}^2)+O(N^{-2\beta})
 \ge c_dN^{-\beta}.
\end{equation}
Since $\OpQ_{\cX_N}(1-p_{d+1}^2)\to d/(d+1)>0$,
\eqref{eq:coordinate-displacement} gives
\[
 \max_j\dist(x_{N,j},y_{N,j})
 \asymp N^{-\beta}\max_j\sqrt{1-p_{d+1}(x_{N,j})^2}
 \asymp N^{-\beta}.
\]
The injectivity of $T_{d+1,\varepsilon_NN^{-\beta}}$ ensures that the
perturbed nodes are distinct.
Since $\OpI(p_{d+1})=0$, for every fixed $s>d/2$ with $s>\beta d$,
\[
 \wce(\OpQ_{\cY_N};H^s)
 \ge\frac{|\OpQ_{\cY_N}(p_{d+1})|}{\norm{p_{d+1}}_{H^s}}
 \gtrsim N^{-\beta}.
\]
When $\beta>1/2$ and $s_{\cX}^*\ge\beta d$, stability gives the reverse
strength bound, so $s_{\cY}^*=\beta d$.
When $0<\beta\le1/2$, every $s>d/2$ satisfies $s>\beta d$, proving
part~(ii).
\end{proof}

\begin{proof}[Proof of Corollary~\ref{cor:accuracy}]
At each fixed $s$, it suffices to assume \eqref{eq:average-s-energy}
if $s<s_c$ or \eqref{eq:logarithmic-packing} if $s=s_c$,
uniformly in $N$; both follow from \eqref{eq:average-energy}.
No energy condition is needed if $s>s_c$.
The definition of $E_s$ and the triangle inequality give
\begin{equation}\label{eq:error-comparison}
 E_s(\cX_N,\delta)
 =\sup_{\cY\in\cP_\delta(\cX_N)}\wce(\OpQ_{\cY};H^s)
 \asymp\wce(\OpQ_{\cX_N};H^s)+\omega_s(\cX_N,\delta).
\end{equation}
Indeed, the triangle inequality gives the upper bound, while
$E_s\ge\wce(\OpQ_{\cX_N};H^s)$ and $\omega_s\le2E_s$ give the lower bound.

The optimal lower bound in
\cite[Theorem~3]{BrauchartSaffSloanWomersley2014} and the reference QMC
property give $\wce(\OpQ_{\cX_N};H^s)\asymp N^{-s/d}$.
Equation~\eqref{eq:error-comparison} and
Lemma~\ref{lem:coordinate-perturbation} prove the lower bound.

For the upper bound when $0<\delta\le1$, the additional subcritical term
satisfies
\[
 N^{-1/2}\delta^{s-d/2}\lesssim N^{-s/d}+\delta.
\]
To see this, split at $\delta=N^{-s/d}$ and use
$N^{-1/2}(N^{-s/d})^{s-d/2-1}
=N^{-(s-1)(s-d/2)/d}\le1$.
At $s=s_c$, the same split gives
\[
 N^{-1/2}\delta\sqrt{1-\ln\delta}
 \le\sqrt{\frac{1+(s/d)\ln N}{N}}
       \max\{N^{-s/d},\delta\}
 \lesssim N^{-s/d}+\delta.
\]
For $s>s_c$, Theorem~\ref{thm:critical-smooth}(ii) applies directly.
For $1<\delta<\pi$, uniform boundedness of point evaluation gives
$E_s\lesssim1\lesssim\delta$; at $\delta=0$, use the reference error.
This proves \eqref{eq:accuracy-balance}.
Finally, applying this estimate at $\min\{\delta_N,\pi/2\}$ shows that
$E_s(\cX_N,\delta_N)=O(N^{-s/d})$ forces $\delta_N\to0$ and then
$\delta_N=O(N^{-s/d})$. The converse is immediate.
\end{proof}

\section{Numerical results}\label{sec:numerics}

We consider equal-weight cubature on $\mathbb S^2$ with three families of numerically computed reference nodes: spherical $t$-designs, maximum-determinant (Fekete) points,
and minimum Coulomb-energy points. The apparent QMC strengths of Fekete and Coulomb-energy points have been estimated as $1.5$ and $2$, respectively, in \cite[Table~1 and Section~8.4]{BrauchartSaffSloanWomersley2014}. Spherical $t$-design families with $N_t\asymp t^2$ have infinite strength. We use $N=81,\allowbreak169,\allowbreak289,\allowbreak441,\allowbreak625,\allowbreak1089,\allowbreak2401,\allowbreak4225$ nodes from each family.

For $\beta\in\{0.4,0.5,0.65,1.8\}$, recall that smaller $\beta$ permits larger perturbations. We perturb the nodes by
\begin{equation}\label{eq:numerical-perturbation}
 y_{N,j}
 =\frac{x_{N,j}+\sin(\varepsilon_NN^{-\beta})e_3}
        {|x_{N,j}+\sin(\varepsilon_NN^{-\beta})e_3|},
 \qquad 1\le j\le N,
\end{equation}
where $e_3=(0,0,1)$ and $\varepsilon_N\in\{-1,1\}$ is chosen so that
\[
 \varepsilon_N\frac1N\sum_{j=1}^N x_{N,j}^{(3)}
 =\left|\frac1N\sum_{j=1}^N x_{N,j}^{(3)}\right|.
\]
This is the perturbation \eqref{eq:signed-perturbation}, with maximal paired displacement of order $N^{-\beta}$. All cubature weights remain $1/N$.

The Sobolev worst-case errors are approximated by the finite sums
\begin{equation}\label{eq:numerical-wce}
 e_{s,L}(\cY_N)^2
 =\sum_{\ell=1}^{L}(1+\ell(\ell+1))^{-s}
   \sum_{k=1}^{2\ell+1}
   \left|\frac1N\sum_{j=1}^N Y_{\ell,k}(y_{N,j})\right|^2,
 \qquad L=4(\sqrt N-1),
\end{equation}
Here $\{Y_{\ell,k}:1\le k\le2\ell+1\}$ is an orthonormal basis of degree-$\ell$ spherical harmonics in $L_2(\mathbb S^2;\sigma_2)$. We take
\[
 s\in\{1.1,1.25,1.5,1.75,2,2.25,2.5,3,3.5,4\}.
\]
For each $s$, we fit the model
\[
 \ln e_{s,L}(\cY_N)\approx a_s-\frac{r(s)}2\ln N,
 \qquad
 e_{s,L}(\cY_N)\approx\e^{a_s}N^{-r(s)/2},
\]
by least squares over the five largest node counts. Here $a_s$ is the
intercept and $-r(s)/2$ is the slope in logarithmic coordinates.
Since $d=2$, the optimal Sobolev rate is $N^{-s/2}$, so $r(s)=s$
corresponds to this rate, while $r(s)<s$ indicates slower observed decay.
For the reference rules and $\beta=0.65,1.8$, increasing the cutoff to
$6(\sqrt N-1)$ at the largest node count changes the fitted exponents
by at most $0.076$, and by less than $0.018$ for $s\ge2$.

\begin{figure}[htbp]
 \centering
 \includegraphics[width=0.75\textwidth]{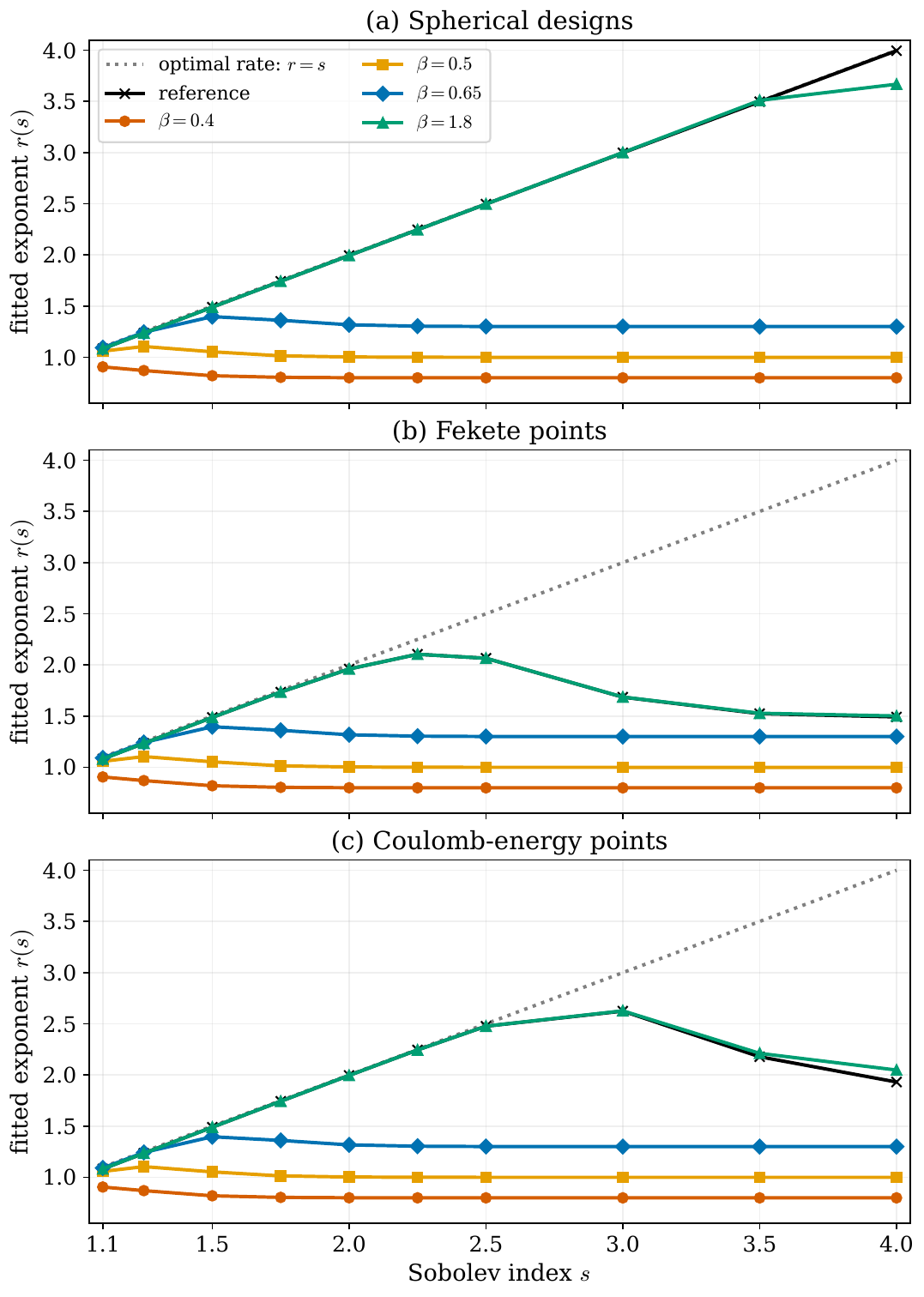}
 \caption{Sobolev error decay under the node perturbation
 \eqref{eq:numerical-perturbation}: (a) spherical designs,
 (b) Fekete points, and (c) Coulomb-energy points.
 The fitted exponent $r(s)$ corresponds to decay $N^{-r(s)/2}$.} 
 \label{fig:theorem-1-1}
\end{figure}

Figure~\ref{fig:theorem-1-1} shows the fitted exponents for the reference
and perturbed configurations. For spherical designs, the reference curve
follows $r(s)=s$ throughout the tested range, as expected from their
infinite QMC strength. At $\beta=0.65$, the exponents for all
three perturbed families level off near $1.3=2\beta$ as $s$ increases,
in agreement with the sharp strength bound in Theorem~\ref{thm:main}(i).
At $\beta=1.8$, the Fekete and Coulomb-energy curves remain close to
their reference curves, consistent with preservation of a reference
strength below $2\beta$. The perturbed spherical designs instead have
a fitted exponent of $3.67$ at $s=4$, close to $2\beta=3.6$.
For $\beta=0.4$ and $0.5$, the fitted exponents lie below $s$ at the
tested smoothness indices and level off near $0.8$ and $1$, respectively.
These rates agree with the loss of admissible QMC rates in
Theorem~\ref{thm:main}(ii).

\FloatBarrier
\bibliographystyle{siamplain}
\bibliography{myref} 

\begin{thebibliography}{10}

\bibitem{AnWu2024Bypassing}
{\sc C.~An and H.-N. Wu}, {\em Bypassing the quadrature exactness assumption of
  hyperinterpolation on the sphere}, J. Complexity, 80 (2024), p.~101789,
  \url{https://doi.org/10.1016/j.jco.2023.101789}.

\bibitem{AtkinsonHan2012}
{\sc K.~Atkinson and W.~Han}, {\em Spherical Harmonics and Approximations on
  the Unit Sphere: An Introduction}, vol.~2044 of Lecture Notes in Mathematics,
  Springer, Berlin, Heidelberg, 2012,
  \url{https://doi.org/10.1007/978-3-642-25983-8}.

\bibitem{BilykDaiMatzke2018}
{\sc D.~Bilyk, F.~Dai, and R.~Matzke}, {\em The {S}tolarsky principle and
  energy optimization on the sphere}, Constr. Approx., 48 (2018), pp.~31--60,
  \url{https://doi.org/10.1007/s00365-017-9412-4}.

\bibitem{BondarenkoRadchenkoViazovska2013}
{\sc A.~V. Bondarenko, D.~Radchenko, and M.~Viazovska}, {\em Optimal asymptotic
  bounds for spherical designs}, Ann. of Math. (2), 178 (2013), pp.~443--452,
  \url{https://doi.org/10.4007/annals.2013.178.2.2}.

\bibitem{BondarenkoRadchenkoViazovska2015}
{\sc A.~V. Bondarenko, D.~Radchenko, and M.~Viazovska}, {\em Well-separated
  spherical designs}, Constr. Approx., 41 (2015), pp.~93--112,
  \url{https://doi.org/10.1007/s00365-014-9238-2}.

\bibitem{BrandoliniEtAl2014Quadrature}
{\sc L.~Brandolini, C.~Choirat, L.~Colzani, G.~Gigante, R.~Seri, and
  G.~Travaglini}, {\em Quadrature rules and distribution of points on
  manifolds}, Ann. Sc. Norm. Super. Pisa Cl. Sci. (5), 13 (2014), pp.~889--923,
  \url{https://doi.org/10.2422/2036-2145.201103_007}.

\bibitem{BrauchartDick2013Sobolev}
{\sc J.~S. Brauchart and J.~Dick}, {\em A characterization of {Sobolev} spaces
  on the sphere and an extension of {Stolarsky}'s invariance principle to
  arbitrary smoothness}, Constr. Approx., 38 (2013), pp.~397--445,
  \url{https://doi.org/10.1007/s00365-013-9217-z}.

\bibitem{BrauchartEtAl2015Covering}
{\sc J.~S. Brauchart, J.~Dick, E.~B. Saff, I.~H. Sloan, Y.~G. Wang, and R.~S.
  Womersley}, {\em Covering of spheres by spherical caps and worst-case error
  for equal weight cubature in {Sobolev} spaces}, J. Math. Anal. Appl., 431
  (2015), pp.~782--811, \url{https://doi.org/10.1016/j.jmaa.2015.05.079}.

\bibitem{BrauchartGrabner2015}
{\sc J.~S. Brauchart and P.~J. Grabner}, {\em Distributing many points on
  spheres: Minimal energy and designs}, J. Complexity, 31 (2015), pp.~293--326,
  \url{https://doi.org/10.1016/j.jco.2015.02.003}.

\bibitem{BrauchartGrabnerSloanWomersley2024Needlets}
{\sc J.~S. Brauchart, P.~J. Grabner, I.~H. Sloan, and R.~S. Womersley}, {\em
  Needlets liberated}, Appl. Comput. Harmon. Anal., 73 (2024), p.~101693,
  \url{https://doi.org/10.1016/j.acha.2024.101693}.

\bibitem{BrauchartHesse2007}
{\sc J.~S. Brauchart and K.~Hesse}, {\em Numerical integration over spheres of
  arbitrary dimension}, Constr. Approx., 25 (2007), pp.~41--71,
  \url{https://doi.org/10.1007/s00365-006-0629-4}.

\bibitem{BrauchartSaffSloanWomersley2014}
{\sc J.~S. Brauchart, E.~B. Saff, I.~H. Sloan, and R.~S. Womersley}, {\em {QMC}
  designs: Optimal order quasi {M}onte {C}arlo integration schemes on the
  sphere}, Math. Comp., 83 (2014), pp.~2821--2851,
  \url{https://doi.org/10.1090/S0025-5718-2014-02839-1}.

\bibitem{DaiXu2013}
{\sc F.~Dai and Y.~Xu}, {\em Approximation Theory and Harmonic Analysis on
  Spheres and Balls}, Springer Monographs in Mathematics, Springer, New York,
  2013, \url{https://doi.org/10.1007/978-1-4614-6660-4}.

\bibitem{DeLaTorreMarzo2024}
{\sc V.~de~la Torre and J.~Marzo}, {\em {QMC} strength for some random
  configurations on the sphere}, in Monte Carlo and Quasi-Monte Carlo Methods,
  A.~Hinrichs, P.~Kritzer, and F.~Pillichshammer, eds., vol.~460 of Springer
  Proceedings in Mathematics \& Statistics, Springer, Cham, 2024, pp.~625--642,
  \url{https://doi.org/10.1007/978-3-031-59762-6_31}.

\bibitem{DelsarteGoethalsSeidel1977}
{\sc P.~Delsarte, J.-M. Goethals, and J.~J. Seidel}, {\em Spherical codes and
  designs}, Geom. Dedicata, 6 (1977), pp.~363--388,
  \url{https://doi.org/10.1007/BF03187604}.

\bibitem{Hesse2006Lower}
{\sc K.~Hesse}, {\em A lower bound for the worst-case cubature error on spheres
  of arbitrary dimension}, Numer. Math., 103 (2006), pp.~413--433,
  \url{https://doi.org/10.1007/s00211-006-0686-x}.

\bibitem{HesseMhaskarSloan2007}
{\sc K.~Hesse, H.~N. Mhaskar, and I.~H. Sloan}, {\em Quadrature in {B}esov
  spaces on the euclidean sphere}, J. Complexity, 23 (2007), pp.~528--552,
  \url{https://doi.org/10.1016/j.jco.2006.10.004}.

\bibitem{HesseSloan2005WCE}
{\sc K.~Hesse and I.~H. Sloan}, {\em Worst-case errors in a {S}obolev space
  setting for cubature over the sphere {${S}^2$}}, Bull. Aust. Math. Soc., 71
  (2005), pp.~81--105, \url{https://doi.org/10.1017/S0004972700038041}.

\bibitem{HesseSloan2006ArbitraryOrder}
{\sc K.~Hesse and I.~H. Sloan}, {\em Cubature over the sphere {${S}^2$} in
  {S}obolev spaces of arbitrary order}, J. Approx. Theory, 141 (2006),
  pp.~118--133, \url{https://doi.org/10.1016/j.jat.2006.01.004}.

\bibitem{SloanWomersley2026}
{\sc I.~H. Sloan and R.~S. Womersley}, {\em {QMC} designs---cubature on the
  sphere without polynomial exactness}, Constr. Approx., 64 (2026),
  pp.~183--208, \url{https://doi.org/10.1007/s00365-026-09760-9}.

\bibitem{Wu2026Hyperinterpolation}
{\sc H.-N. Wu}, {\em Hyperinterpolation beyond exact cubature: A spectral
  multiplier approach}, 2026, \url{https://arxiv.org/abs/2605.17739}.

\bibitem{WuZhuang2026NestedQMC}
{\sc H.-N. Wu and X.~Zhuang}, {\em Nested {QMC} designs on spheres}, 2026,
  \url{https://arxiv.org/abs/2609.15767}.

\bibitem{Zhou2020Neighborhood}
{\sc Y.~Zhou}, {\em Quadrature rules with neighborhood of spherical designs on
  the two-sphere}, Appl. Math. Comput., 367 (2020), p.~124769,
  \url{https://doi.org/10.1016/j.amc.2019.124769}.

\bibitem{ZhouChen2018}
{\sc Y.~Zhou and X.~Chen}, {\em Spherical {$t_\varepsilon$}-designs for
  approximations on the sphere}, Math. Comp., 87 (2018), pp.~2831--2855,
  \url{https://doi.org/10.1090/mcom/3306}.

\end{thebibliography}

\end{document}